\documentclass[12pt]{article}
\usepackage{amssymb}

\usepackage{amsthm}
\usepackage{amsmath}
\usepackage[T1]{fontenc}
\usepackage[utf8]{inputenc}
\usepackage[english]{babel}
\usepackage[pdftex]{graphicx}          

\usepackage{indentfirst}              
\usepackage{makeidx}                    
\usepackage[nottoc]{tocbibind}         
\usepackage{courier}                   
\usepackage{type1cm}                    
\usepackage{listings}                
\usepackage{titletoc}
\usepackage[all,cmtip]{xy}

\newcommand{\HH}{\mathbb{H}}
\newcommand{\OO}{\mathbb{O}}
\newcommand{\Z}{\mathbb{Z}}

\newcommand{\Cay}{{\rm Cay\,}}
\newcommand{\Ca}{{\rm Cay}}
\newcommand{\Reg}{{\rm Reg\,}}

\usepackage{fancyhdr}
\newtheorem{defi}{\bf Definition}[section]
\newtheorem{theo}[defi]{\bf Theorem}
\newtheorem{lemma}[defi]{\bf Lemma}
\newtheorem{coro}[defi]{\bf Corollary}
\newtheorem{pro}[defi]{\bf Proposition}

\def\bes{\begin{eqnarray*}}
\def\ees{\end{eqnarray*}}
\def\bee{\begin{eqnarray}}
\def\eee{\end{eqnarray}}
\def\a{\alpha}
\def\b{\beta}

\def\f{\varphi}
\def\la{\langle}
\def\ra{\rangle}

\newenvironment{dem}{\noindent
\textnormal{\textbf{Proof:}}}{{\hfill$\Box$}}

\theoremstyle{remark}
\newtheorem{remark}[defi]{\bf Remark}

\newcommand{\algA}{\mathcal{A}}
\newcommand{\algB}{\mathcal{B}}
\newcommand{\algC}{\mathcal{C}}

\title{On a problem by Nathan Jacobson}

\author{
Victor H.\,L\'opez Sol\'{\i}s\\{\small   Departamento Acad\'emico de Matem\'atica}\\
{\small Facultad de Ciencias}\\
{\small Universidad Nacional Santiago Ant\'unez de Mayolo}\\
{\small Huaraz, Per\'u}\\
{\small vlopezs@unasam.edu.pe}
\and
Ivan P.\,Shestakov \\{\small Instituto de Matem\'atica e Estat\'{\i}stica}\\
{\small Universidade de S\~ao Paulo},\\ {\small S\~ao Paulo. Brasil}\\
{\small and}\\
{\small Sobolev Institute of Mathematics,}\\ {\small Novosibirsk, Russia}\\
{\small ivan.shestakov@gmail.com}
}
\date{\quad}

\begin{document}
\maketitle\vspace{-1.5cm}

\begin{abstract}
We prove a coordinatization theorem for unital alternative algebras containing $2\times 2$ matrix algebra with the same identity element 1. This solves an old problem announced by Nathan Jacobson \cite{J1} on the description of alternative algebras containing a generalized quaternion algebra $\mathbb{H}$ with the same 1, for the case when the algebra $\HH$ is split. In particular, this is the case when the basic field is finite or algebraically closed.
\end{abstract}

{\parindent= 4em \small  \sl Keywords: Alternative algebras, Kronecker factorization theorem, quaternion algebra, Cayley bimodule, Pl\"ucker relations}

\section{Introduction.}

The classical Wedderburn Coordinatization Theorem says that if a unital associative algebra $A$ contains a matrix algebra $M_n(F)$  with the same identity element then it is itself a matrix algebra, $A\cong M_n(D)$, ``coordinated'' by $D$. Generalizations and analogues of this theorem were proved for various classes of algebras and superalgebras \cite{J1,K,LDS1,LDS2,  MSZ,CE,Mac,PS,Sc}. The common content of all these results is that if an algebra (or superalgebra) contains a certain subalgebra (matrix algebra, octonions, Albert algebra) with the same unit then the algebra itself has the same structure,  but not over the basic field rather over a certain algebra that ``coordinatizes'' it.  The Coordinatization Theorems play important role in structure theories, especially in classification theorems, and also in the representation theory, since quite often an algebra $A$ coordinated by $D$  is Morita equivalent to $D$, though they could belong to different classes (for instance, Jordan algebras are coordinated by associative and alternative algebras).

\medskip  

In this paper we consider alternative algebras. Recall that an algebra $\algA$ is called {\em alternative} if it  satisfies the following identities:
\begin{equation}\label{e1}
x^2 y = x(xy), ~~ (xy)y = xy^2. 
\end{equation}
for all $x,y \in\algA$.
All associative algebras are clearly alternative.
A classical  example of a non-associative alternative algebra is the Cayley  (or generalized octonion) algebra $\mathbb{O}$ (see \cite{J2, Sc, KS, ZSSS}). 
Kaplansky \cite{K} proved an analogue  of Wedderburn's theorem for alternative algebras containing the Cayley algebra. He showed
that if $\algA$ is an alternative algebra with identity element $1$ which contains a subalgebra $\mathcal{B}$
isomorphic to a Cayley algebra and if $1$ is contained in $\mathcal{B}$, then $\algA$ is isomorphic to the Kronecker product $\mathcal{B}\otimes T$, where $T$ is the center
of $\algA$. Jacobson gave a new proof of Kaplansky's result, using his classification of irreducible alternative bimodules, and in addition proved an analogue of this theorem for Jordan algebras \cite{J1}, where the role of Cayley algebra is played by the Albert algebra, the exceptional simple Jordan algebra of dimension $27$. These results have important applications in the theory of representations of alternative and Jordan algebras \cite {J, J2}.
\smallskip

The Wedderburn coordinatization theorem in the case $n\geq 3$ admits a generalization for  alternative algebras, since every alternative algebra $\algA$ which contains a subalgebra $M_n(F)$ $(n\geq 3)$ with the same identity element is associative (see \cite[Corollary $11$, Chapter $2$]{Sc}). The result is not true for $n=2$, the split Cayley algebra and its 6-dimensional subalgebra are counterexamples. The problem of description of alternative algebras containing $M_2(F)$ or, more generally,  a generalized quaternion algebra $\mathbb{H}$ with the same identity element was posed by Jacobson \cite{J1}. In this paper, we solve this problem for the split case $\mathbb{H}\cong M_2(F)$, without any restriction on the dimension and characteristic of the base field $F$.
\smallskip

Our $M_2(F)$-coordinatization involves two ingredients: an alternative $M_2(F)$-algebra $\algA$ is ``coordinated'' by an associative algebra $\algB$ and by a commutative $\algB$-bimodule $V$ (that is, $V$ is annihilated by any commutator of elements of $\algB$), on which a skew-symmetric mapping is defined with values in the center of $\algB$, satisfying Pl\"ucker relations. More exactly,  $\algA=M_2(\algB)\oplus V^2$, with a properly defined multiplication. The details are given in the Main Theorem \ref{th_main}.
\smallskip

The paper is organized as follows. In Sections 1-2 we give definitions and some known results on alternative algebras ans bimodules. In Section 3 we prove that a unital alternative algebra $\algA$ containing the generalized quaternion algebra $\HH$ with the same unit admits a $Z_2$-grading $\algA=\algA_a\oplus\algA_c$ with associative $0$-component $\algA_a$.  In the next section we determine multiplication in the 1-component $\algA_c$. In Section 5 we prove The Main Theorem on $M_2(F)$-coordinatization of alternative algebras. Section 6-7 are devoted to examples and open questions.
\smallskip 

Throughout this paper the ground field F is of arbitrary characteristic.
\section{Definitions and known results}

Let $\algA$ be an arbitrary algebra. Denote by $(x,y,z)=(xy)z-x(yz)$ the \textit{associator} of the elements $x,y,z\in\algA$ and by $[x,y]=xy-yx$ the \textit{commutator} of the elements $x,y\in\algA$. For subsets $\mathcal{B}, \mathcal{C}, \mathcal{D}$ of $\algA$, we denote by  $(\mathcal{B},\mathcal{C},\mathcal{D})$  the \textit{ associator space} generated by all the associators $(b,c,d)$, $b\in\mathcal{B}$, $c\in\mathcal{C}$, $d\in\mathcal{D}$.
The \textit{associative center} $N(\algA)$, the \textit{commutative center} $K(\algA)$ and the \textit{center} $Z(A)$ are respectively defined as follows:
\begin{gather}
N(\algA)=\{a\in\algA|(a,\algA,\algA)=(\algA,a,\algA)=(\algA,\algA,a)=0\}, \notag
\\[4pt]
K(\algA)=\{a\in\algA|[a,\algA]=0\},\notag
\\[4pt]
Z(\algA)=N(\algA)\cap K(\algA).\notag
\end{gather}

In term of  associators,  identities $(\ref{e1})$ defining  alternative algebras can be written as
$$(x,x,y)=0,~~(x,y,y)=0.$$
The first of them is called the \textit{left alternative identity} and the  second one, the \textit{right alternative identity}.

Linearizing the left and right alternative identities, we obtain
$$(x,z,y)+(z,x,y)=0,\ \,\,\,(x,y,z)+(x,z,y)=0,$$ 
which show that in an alternative algebra the associator is an antisymmetric function of its arguments. Also, these identities can be written as
\begin{equation}\label{e2}
(x\circ z)y-x(zy)-z(xy)=0,~~(xy)z+(xz)y-x(y\circ z)=0,
\end{equation}
where $a\circ b=ab+ba$.

Throughout the article we will make use of some identities that are valid in any alternative algebra and will be mentioned at the time be required.

\subsection{Alternative bimodules}

Let be $\algA$ an alternative algebra over $F$ and $V$ a bimodule over $\algA$, this is, $V$ is a vector space over $F$ equipped with the applications $\algA\otimes V\longrightarrow V$, $a\otimes v\longmapsto av$, $V\otimes \algA\longrightarrow V$, 
$v\otimes a\longmapsto va$, $a\in\algA$, $v\in V$. Define on the vector space $E=\algA\oplus V$ a binary operation  $\cdot:E\times E\longrightarrow E$ by 
$$(a+v)\cdot(b+w)=ab+ (av+wb),$$
where $a,b\in\algA$, $v,w\in V$. Then $E$ with the operation (product) $\cdot$ becomes an algebra, the \textit{\textbf{split null extension}} of $\algA$ by bimodule $V$, where $\algA$ is a subalgebra and $V$ is an ideal such that $V^{2}=0$.  Now, $V$ is called an \textit{\textbf{alternative bimodule}} over $\algA$ if
$E$ is an alternative algebra with respect to $\cdot$. 

Due to identities (\ref{e1}),  a bimodule $V$ over $\algA$ is an alternative bimodule if and only if the following relationships are satisfied:
\begin{gather}
(a,a,v)=0,~~(a,v,b)+(v,a,b)=0, \notag
\\[5pt]
(v,b,b)=0,~~(a,v,b)+(a,b,v)=0, \notag
\end{gather}
for all $a,b\in\algA$, $v\in V$.

Let $\algA$ be a composition algebral (see \cite{J2, KS, Sc, ZSSS}). Recall that $\algA$ is a unital alternative algebra, it has an involution $a\mapsto a^*$ such that the {\em\textbf {trace}} $t(a)=a+a^*$ and {\em\textbf {norm}} $n(a)=aa^*$ lie in $F$.

An alternative bimodule $V$ over a composition algebra $\algA$ is called a \textit{\textbf{Cayley bimodule}} if it satisfies the relation
\begin{equation}\label{e3}
av=v{a^*},
\end{equation}
where $a\in\algA$, $v\in V$, e $a\rightarrow {a^*}$ is the canonical involution in $\algA$.

\smallskip
Typical examples of composition algebras are the algebras of (generalized) quaternions $\HH$ and octonions $\OO$ with symplectic involutions.
 Recall that $\OO=\HH\oplus v\HH$, with the product defined by
 \bee\label{CD_product}
{{a\cdot b=ab,\ a\cdot vb=v(a^*b),\ vb\cdot a=v(ab), \ va\cdot vb= (ba^*)v^2,}}
\eee
where $a,b\in\HH,\ 0\neq v^2\in F, \ a\mapsto a^*$ is the symplectic involution in $\HH$.

The subspace {{$v\HH\subset\OO$}} is invariant under multiplication by elements of $\HH$ and it gives an example of a {{\em Cayley bimodule}} over $\HH$.
If  {{$\HH$ is a division algebra then $v\HH$ is irreducible, otherwise $\HH\cong M_2(F)$ and  
\[
v\HH=\langle ve_{\tiny{22}},-ve_{\footnotesize{12}}\rangle\oplus\langle -ve_{\small{21}},ve_{11}\rangle,
 \] 
 where $M_2(F)$-bimodules $\langle ve_{\tiny{22}},\,-ve_{\footnotesize{12}}\rangle$ and $\langle -ve_{\small{21}},ve_{11}\rangle$ are both isomorphic to the 2-dimensional Cayley bimodule
$\Ca=F\cdot m_{1}+F\cdot m_{2}$, with the action of $M_{2}(F)$ given by
\bee\label{id_cay}
e_{ij}\cdot m_{k}=\delta_{ik}m_{j},\ \ \ \  
m\cdot a={a^*}\cdot m, 
\eee
where $a\in M_{2}(F),\, m\in \Ca,\ i,j,k\in\{1,2\}$ and $a\mapsto {a^*}$ is the symplectic involution in $M_{2}(F)$.

   We will denote the Cayley bimodule $v\HH$ for division $\HH$ as {{$\Cay \HH$}}, and the regular (associative) $\HH$-bimodule
by {{$\Reg$}}.

\section{$\Z_2$-grading  $\algA=\algA_a+\algA_c$.}

The statement  of the following result follows from  \cite[Lemma 11]{Sh1} and its proof.

\begin{pro}\label{pr1}{\rm \cite[Lemma 11]{Sh1}}.
Let $\algA$ be a unitary alternative algebra over the field $F$ which contains a composition subalgebra $\mathcal{C}$ with the same identity element. Suppose that a subspace $V$ of $\algA$ is $\mathcal{C}-$invariant and satisfies $(\ref{e3})$. Then, the following identities are valid for any $a,b\in\mathcal{C}$, $r\in\algA$, $u,v\in V$;
\begin{gather}
(ab)v=b(av), ~v(ab)=(vb)a, \label{e4}
\\[5pt]
a(ur)=u({a^*}r), \label{e5}
\\[5pt]
a(uv)=u(va), ~(uv)a=(au)v, \label{e6}
\\[5pt]
(u,v,a)=[uv,a]  \label{e7}
\end{gather}
\end{pro}

It is important to know the structure of unitary alternative $\mathbb{H}-$bimodules. Their structure is given by the following result:

\begin{pro}\label{pr2}
{\rm {\cite[Lemma 12]{Sh1}}}. Every unitary alternative $\mathbb{H}-$bimodule $V$ is completely reducible and admits decomposition $V=V_{a}\oplus V_{c}$, where $V_{a}$ is an associative $\mathbb{H}-$bimodule and $V_{c}$ is a Cayley bimodule over $\mathbb{H}$; furthermore, the subbimodule $V_{c}$ coincides with the associator subspace $(V,\mathbb{H},\HH)$.
Every irreducible component of the subbimodule $V_{a}$ is isomorphic to the regular $\HH$-bimodule $\Reg$, and every irreducible component of the subbimodule $V_{c}$ is isomorphic to $\Cay \HH$ if $\HH$ is a division algebra and to $\Ca$ if $\HH\cong M_2(F)$.
\end{pro}

Let $\algA$ be an alternative algebra such that $\algA$ contains $\mathbb{H}$ with the same identity element, so we can consider $\algA$ as a unitary alternative $\mathbb{H}-$bimodule. Then, by Proposition $\ref{pr2}$, $\algA$ is completely reducible and admits the decomposition 
$$
\algA=\algA_{a}\oplus\algA_{c},
$$
where $\algA_{a}$ is a unitary associative  $\mathbb{H}-$bimodule and $\algA_{c}$ is a unitary Cayley $\mathbb{H}-$bimodule.

\smallskip
Denote $Z_{a}=\{u\in\algA_{a}|[u,\mathbb{H}]=0\}$. Since $\algA_{a}$ is  isomorphic to a direct sum of bimodules $\Reg$, we have  $\algA_{a}=\sum_{i}\oplus\Reg_{i},\ \Reg_i\cong\Reg$ for all $i$. This implies that $\algA_{a}$ contains a set of elements $\{u_{i}\}$ (the images of $1$ under the isomorphisms with $\mbox{Reg}$) such that $\mbox{Reg}_{i}=u_{i}\mathbb{H}$ with $u_{i}\in Z_{a}$, and each element of $\algA_{a}$ can be written in only one way in the form $\sum u_{i}a_{i}$, $a_{i}\in \mathbb{H}$. Now, of course, $Z_{a}\neq 0$ and $\algA_{a}=Z_{a}\mathbb{H}$.

Also, by Proposition $\ref{pr2}$, the bimodule $\algA_{c}$ coincides with $(\algA,\mathbb{H},\mathbb{H})$ and is completely reducible; this is, $\algA_{c}=\sum_{j}\oplus\widetilde{\mbox{Cay}}_{j}$, where $\widetilde{\mbox{Cay}}$ is equal to $\mbox{Cay}\,\mathbb{H}$ or to $\mbox{Cay}$. Therefore,
$$
\algA=(\sum\oplus\mbox{Reg}_{i})\oplus(\sum\oplus\widetilde{\mbox{Cay}}_{j}).
$$

The statements and demonstrations of Lemmas $\ref{l1}$ and $\ref{l2}$ are similar to Lemmas $3.1$ and $3.2$ of $\cite{LDS2}$ given there for superbimodules over the superalgebra $B(4,2)=\HH+\Ca$.

\begin{lemma}\label{l1}
Let $\algA=\algA_{a}\oplus\algA_{c}$ be the decomposition of $\algA$ from above. Then for any $m,n\in\algA_{c}$, $a\in \mathbb{H}$,
\begin{equation}\label{e8}
(mn)a=(am)n, ~~~ a(mn)=m(na) 
\end{equation}
and for any $u\in \algA_{a}$, $m\in \algA_{c}$, $a,b\in \mathbb{H}$,
\bee
(um)a&=&(u{a^*})m \label{e9}\\
a(mu)&=&m({a^*}u), \label{e10}
\\
((um)a)b&=&(um)(ba), \label{e11}
\\
b(a(mu))&=&(ab)(mu), \label{e12}
\\
(um,a,b)&=&(um)[b,a], \label{e13}
\\
(b,a,mu)&=&[b,a](mu). \label{e14}
\eee
\end{lemma}
\begin{dem}
First, let us consider $m,n\in\algA_{c}$, $a\in \mathbb{H}$. By $(\ref{e3})$, $(mn)a-(am)n=(mn)a-
(m{a^*})n = (m, n, a)-(m,{a^*},n)+m(na-{a^*}n)=(m,n,a)+(m,a,n)=0$. Analogously $a(mn)-m(na)=0$.
This proves $(\ref{e8})$.

Now let $u\in \algA_{a}$, $m\in \algA_{c}$, $a,b\in\mathbb{H}$. Then 
$$
(um)a-(u{a^*})m=(u,m,a)-(u,{a^*},m)+u(ma-{a^*}m)=0,
$$ 
and similarly $a(mu)-m({a^*}u)=0$, which proves $(\ref{e9})$ and $(\ref{e10})$.
In addition by $(\ref{e9})$
$$(um)a.b=(u{a^*}.m)b=(u{a^*}.{b^*})m=(u.(ba)^*)m=(um)(ba),$$
which proves $(\ref{e11})$. Similarly, by $(\ref{e10})$, we get $(\ref{e12})$. Finally, by $(\ref{e11})$ and $(\ref{e12})$ we have
\bes
(um,a,b)&=&((um)a)b-(um)(ab)=(um)(ba)-(um)(ab)=(um)[b,a]\\
(b,a,mu)&=&(ba)(mu)-b(a(mu))=(ba)(mu)-(ab)(mu)=[b,a](mu),
\ees
which proves $(\ref{e13})$ and $(\ref{e14})$.
\end{dem}

\begin{lemma}\label{l2}
The products $\algA_{a}\algA_{a}$, $\algA_{a}\algA_{c}$, $\algA_{c}\algA_{a}$ and $\algA_{c}\algA_{c}$ 
are $\mathbb{H}-$invariants subspaces. Moreover $\algA_{a}\algA_{c}+\algA_{c}\algA_{a}\subseteq \algA_{c}$ and $\algA_{c}\algA_{c}\subseteq \algA_{a}$.
\end{lemma}
\begin{dem}
Since $\algA_{a}$ and $\algA_{c}$ are $\mathbb{H}-$invariante, in order  to prove  the first part of the Lemma it suffices to show that the product of any $\mathbb{H}-$invariants subspaces $U$ and $W$ is again  $\mathbb{H}-$invariant.

We have by the linearized identity of the right alternativity $(\ref{e2})$
$$
(UW)\mathbb{H}\subseteq U(W\circ \mathbb{H})+(U\mathbb{H})W\subseteq UW,
$$ 
and similarly $\mathbb{H}(UW)\subseteq UW.$

Now, let us demonstrate that $\algA_{a}\algA_{c}+\algA_{c}\algA_{a}\subseteq \algA_{c}$. Recall that by Proposition $\ref{pr2}$,
$\algA_{c}=(\algA,\mathbb{H},\mathbb{H})$. Choose $a,b\in \mathbb{H}$ such that $0\neq [a, b]^2\in F$, then by $(\ref{e14})$
$$
\algA_{c}\algA_{a} = [a,b]^2 (\algA_{c}\algA_{a} )\subseteq [a,b](\algA_{c}\algA_{a})= (a,b,\algA_{c}\algA_{a})\subseteq (\mathbb{H},\mathbb{H},\algA)=\algA_{c},
$$
and similarly $\algA_{c}\algA_{a}\subseteq \algA_{c}$.
Finally, for any $m,n\in\algA_{c}$ and $a\in \mathbb{H}$, we have by $(\ref{e8})$ and $(\ref{e4})$
$$((mn)a)b = ((am)n)b = (b(am))n = ((ab)m)n = (mn)(ab),$$
which proves $\algA_{c}\algA_{c}\subseteq \algA_{a}.$
\end{dem}

\begin{lemma}\label{l3}
$\algA_{a}$ is an associative subalgebra of $\algA$.
\end{lemma}
\begin{dem}
 Recall the following identities valid in every alternative algebra (see \cite{ZSSS, Sh1}).
\bee
(xy)(zx)&=&x(yz)x,\label{Moufangm}\\ \ 
[x,yz]&=&[x,y]z+y[x,z]-3(x,y,z), \label{e15}
\\
(xy,z,t)&=&x(y,z,t)+(x,z,t)y-(x,y,[z,t]), \label{e16}
\\
2[(x,y,z),t]&=&([x,y],z,t)+([y,z],x,t)+([z,x],y,t),\label{e17} 
\\ \ 
[x,y](x,y,z)&=&(x,y,(x,y,z))=-(x,y,z)[x,y],\label{e18} 
\\
((z,w,t),x,y)& =&((z,x,y),w,t)+(z,(w,x,y),t) \label{e19}\\
&+&(z,w,(t,x,y))-[w,(z,t,[x,y])]+([z,t],w,[x,y]).
\eee
Let us fix arbitrary elements $u,v,w\in Z_{a}$ and $a,b,c\in \mathbb{H}$. Then by $(\ref{e17})$
$$([a,b],u,v)=2[(a,b,u),v]-([b,u],a,v)-([u,a],b,v)=0.$$
So by $(\ref{e16})$
$$(uv,a,b)=u(v,a,b)+(u,a,b)v-(u,v,[a,b])=-([a,b],u,v)=0,$$
which implies $(Z_{a}Z_{a},\mathbb{H},\mathbb{H})=0$. 
By  linearization of $(\ref{e18})$, choosing $a,b\in \mathbb{H}$ such that $[a,b]^{2}=\alpha\in F$, $\alpha\neq 0$, we have for any $x\in\algA_{a}$
\begin{equation*}
\begin{split}
[a,b](u,x,c)& =-[u,b](a,x,c)-[a,x](u,b,c)-[u,x](a,b,c)+(a,b,(u,x,c))\\
& \quad +(u,b,(a,x,c))+(a,x,(u,b,c))+(u,x,(a,b,c))\\
& =(a,b,(u,x,c))\\
&\overset{(\ref{e19})}{=}((u,a,b),x,c)+(u,(x,a,b),c)+(u,x,(c,a,b))\\
&\quad -[x,(u,c,[a,b])]+([u,c],x,[a,b])=0.
\end{split}
\end{equation*}
So $\alpha(u,x,c)=[a,b]^{2}(u,x,c)=[a,b]([a,b](u,x,c))=0$, thus $(u,x,c)=0$, which implies
$(Z_{a},\algA_{a},\mathbb{H})=0$. In particular, $(Z_{a},Z_{a},\mathbb{H})=0$. Then, by $(\ref{e15})$ 
$$
[a,uv]=[a,u]v+u[a,v]-3(a,u,v)=0,
$$
and so $[\mathbb{H},Z_{a}Z_{a}]=0$. Therefore
$Z_{a}Z_{a}\subseteq Z_{a}$.
\smallskip

By linearization of $(\ref{e18})$, we have
\begin{equation*}
\begin{split}
[a,b](u,v,w)& =-[a,v](u,b,w)-[u,b](a,v,w)-[u,v](a,b,w)+(a,b,(u,v,w))\\
& \quad +(a,v,(u,b,w))+(u,b,(a,v,w))+(u,v,(a,b,w))=0
\end{split}
\end{equation*}
Choose again $a,b\in \mathbb{H}$ such that $[a,b]^{2}=\alpha\in F$, $\alpha\neq 0$. Then
$$\alpha(u,v,w)=[a,b]^{2}(u,v,w)=[a,b]([a,b](u,v,w))=0,$$ and so $(u,v,w)=0$. Thus $Z_{a}$ is an associative algebra.

Consequently, by linearization of the central Moufang identity $(\ref{Moufangm})$ and using the fact that $\algA_{a}$ is an associative $\mathbb{H}-$bimodule, we have
\begin{equation*}
\begin{split}
(ua)(vb)& =-(ba)(vu)+(u(av))b+(b(av))u\\
            &= -(ba)(vu)+(u(va))b+((ba)v)u\\
			&= -(ba)(vu)+((uv)a)b+(ba)(vu)\\
			&= (uv)(ab).
\end{split}
\end{equation*}
Therefore
$\algA_{a}\algA_{a}\subseteq\algA_{a}$, that is, $\algA_{a}$ is a subalgebra of $\algA$.

Remembering that
$(Z_{a},\algA_{a},\mathbb{H})=0$, then for all $x,y\in\algA_{a}$ we have
\begin{equation*}
\begin{split}
(ua,x,y)&\overset{(\ref{e16})}{=}u(a,x,y)+(u,x,y)a-(u,a,[x,y])\\
                            &\,= u(a,x,y)+(u,x,y)a.
\end{split}
\end{equation*}														
Thus, using the last equality several times and the fact that $Z_{a}$ is associative, we have
\begin{equation*}
\begin{split}
(ua,vb,wc)&=u(a,vb,wc)+(u,vb,wc)a\\
         &=-u(v(b,a,wc)+(v,a,wc)b)-(v(b,u,wc)+(v,u,wc)b)a\\
		 &=-((w(c,v,u)+(w,v,u)c)b)a=0.
\end{split}
\end{equation*}
Therefore, $\algA_{a}$ is an associative subalgebra of $\algA$.

\end{dem}

It follows immediately from Lemmas $\ref{l2}$ and $\ref{l3}$, the following result.

\begin{coro}\label{c1}
$\algA=\algA_{a}\oplus\algA_{c}$ is a $\mathbb{Z}_2-$graded algebra, where $\algA_{a}$ is the even part and $\algA_{c}$ is the odd part of the $\mathbb{Z}_2-$grading of $\algA$.
\end{coro}

In what follows we will use in a permanent way the fact that $\algA$ is an $\mathbb{Z}_{2}-$graded alternative algebra. Thus, we have $(\algA_{c},\mathbb{H},\algA_{c})\subseteq(\algA_{c},\algA_{a},\algA_{c})\subseteq\algA_{a}$,  and $[Z_{a},\algA_{c}]\subseteq\algA_{c}$.

\begin{lemma}\label{l4}
$[Z_a,\algA_c]=(Z_{a},\algA,\algA)=0$.
\end{lemma}
\begin{dem}
Let us fix arbitrary elements $u,v,w\in Z_{a}$, $m,n\in\algA_{c}$ and $a,b,c\in \mathbb{H}$. In the proof of the previous Lemma we have shown that $(Z_{a},\algA_{a},\mathbb{H})=0$. So, let us generalize the previous equality and show first
\begin{equation}\label{e20}
(Z_{a},\algA,\mathbb{H})=0.
\end{equation}
By the fact that $\algA_{c}$ is a Cayley $\mathbb{H}-$bimodule, we have
$$(a,u,m)=(au)m-a(um)=(au)m-(um){a^*}=(au)m-(ua)m=[a,u]m=0,$$
which proves $(\mathbb{H},Z_{a},\algA_{c})=0$. Thus,
 $$(Z_{a},\algA,\mathbb{H})\subseteq (Z_{a},\algA_{a},\mathbb{H})+(Z_{a},\algA_{c},\mathbb{H})=0,$$
which proves $(\ref{e20})$. 

In addition, consider the identity
\begin{equation}\label{e21}
([x,y],y,z)=[y,(x,y,z)]
\end{equation}
which is valid in every alternative algebra. Using its linearization, we obtain
$$([u,m],a,b)=-([u,a],m,b)+[m,(u,a,b)]+[a,(u,m,b)]=0.$$
Thus $([Z_{a},\algA_{c}],\mathbb{H},\mathbb{H})=0$, that is, $[Z_{a},\algA_{c}]\subseteq\algA_{a}$. Therefore
$$[Z_{a},\algA_{c}]\subseteq\algA_{a}\cap\algA_{c}=0,$$
which implies $[Z_{a},\algA_{c}]=0$.

By  linearization of $(\ref{e18})$, by $(\ref{e19})$, choosing $a,b\in \mathbb{H}$ such that $0\neq [a,b]^{2}=\alpha\in F$,  we have
\begin{equation*}
\begin{split}
[a,b](u,m,n)&=-[u,b](a,m,n)-[a,m](u,b,n)-[u,m](a,b,n)+(a,b,(u,m,n))\\
                    &\quad +(u,b,(a,m,n))+(a,m,(u,b,n))+(u,m,(a,b,n))\\
										&=	(u,m,(a,b,n))\\
										&=((u,m,a),b,n)-((u,b,n),m,a)-(u,(m,b,n),a)\\
										&\quad +[m,(u,a,[b,n])]-([u,a],m,[b,n])=0,
\end{split}
\end{equation*}
hence $\alpha(u,m,n)=[a,b]^{2}(u,m,n)=[a,b]([a,b](u,m,n))=0$ and $(u,m,n)=0$; therefore, $(Z_{a},\algA_{c},\algA_{c})=0$.
Also
\begin{equation*}
\begin{split}
[a,b](u,n,v)&=-[u,b](a,n,v)-[a,n](u,b,v)-[u,n](a,b,v)+(a,b,(u,n,v))\\
                    &\quad +(u,b,(a,n,v))+(a,n,(u,b,v))+(u,n,(a,b,v))\\
										&=	(a,b,(u,n,v))\\
										&=-(u,(a,n,v),b)-(u,a,(b,n,v))+((u,a,b),n,v)\\
										&\quad +[a,(u,b,[n,v])]-([u,b],a,[n,v])=0,
\end{split}
\end{equation*}
so $\alpha(u,n,v)=[a,b]^{2}(u,n,v)=[a,b]([a,b](u,n,v))=0$ and $(u,n,v)=0$; thus, $(Z_{a},\algA_{c},Z_{a})=0$.
Then by $(\ref{e16})$ and $(\ref{e20})$
$$(ua,v,m)=u(a,v,m)+(u,v,m)a-(u,a,[v,m])=0;$$
so $(Z_{a}\mathbb{H},Z_{a},\algA_{c})=0$. Therefore $(\algA_{a},Z_{a},\algA_{c})=0$, and we have
$$
(Z_{a},\algA,\algA)\subseteq (Z_{a},\algA_{a},\algA_{a})+(Z_{a},\algA_{c},\algA_{a})+(Z_{a},\algA_{c},\algA_{c})=0,
$$
so  $Z_{a}\subseteq N(\algA)$.

\end{dem}

\begin{coro}\label{c2}
$\algA_a=Z_{a}\otimes_{F} \mathbb{H}$.
\end{coro}
\begin{dem}
As $\algA_{a}=\sum \oplus u_i\mathbb{H}$, every element of $\algA_a$ can be written uniquely in the form $\sum u_i a_i$ with $a_i\in\mathbb{H}$. We know that $\algA_a$ is associative. On the other hand, let $x=\sum u_i a_i\in Z_a$ be then $ax=xa$ for all $a\in\mathbb{H}$. Therefore, by $[Z_a,\mathbb{H}]=0$ we have
$$\sum u_i aa_i=\sum u_i a_ia;$$
so, $aa_i=a_ia$. But as $\mathbb{H}$ is central, we have $a_i=\alpha_i 1$, $\alpha_i\in F$. Then $Z_a=\sum Fu_i$ and $\algA_a=Z_{a}\otimes_{F} \mathbb{H}$.

\end{dem}

\begin{lemma}\label{l5}
$[Z_{a},Z_{a}]\algA_{c}=\algA_{c}[Z_{a},Z_{a}]=0.$
\end{lemma}
\begin{dem}
In the proof of Lemma $\ref{l4}$ we have obtained $[Z_{a},\algA_{c}]=0$. Thus, by $(\ref{e15})$ and again by Lemma $\ref{l4}$
$$[Z_{a},Z_{a}]\algA_{c}\subseteq [Z_{a},Z_{a}\algA_{c}]-Z_{a}[Z_{a},\algA_{c}]+3(Z_{a},Z_{a},\algA_{c})=0$$
and similarly $\algA_{c}[Z_{a},Z_{a}]=0.$
\end{dem}

\begin{remark}\label{re2}
Note that in general $Z_a$ is not commutative. For example, if $\algA=M_{n}(\HH)$ then $Z_a\cong M_n(F)$.
If $\algA$ is prime and nonassociative then by \cite[Corollary to Theorem 8.11]{ZSSS} $N(\algA)=Z(\algA)$, hence $Z_a\subseteq Z(\algA)$ is commutative. In fact, in this case $A$ is a Cayley-Dickson ring (see \cite{ZSSS}).
\end{remark}

\section{Multiplication in $\algA_c$.}

In the previous section we described, in particular, the structure of the associative part $\algA_a$. This section is devoted to description of the multiplication in the Cayley part $\algA_c$. 
Here and below we will assume that the quaternion algebra $\HH$ is split, that is, $\HH\cong M_2(F)$. 
 
 \smallskip
 
We have already mentioned that the Cayley $\HH$-bimodule $\algA_c$ is completely reducible and is a direct sum of bimodules isomorphic to the Cayley bimodule $\Ca= F\cdot m_1+F\cdot m_2$ from \eqref{id_cay}. Denote by $V(1)$ and $V(2)$ the subspaces of $\algA_c$ spanned by the elements of type $m_1$ and $m_2$, respectively; then the mappings 
\bes
&\pi_{12}:V(1)\rightarrow V(2),\ v\mapsto e_{12}\cdot v, &\\
 &\pi_{21}:V(2)\rightarrow V(1),\ v\mapsto  e_{21}\cdot v&
 \ees
are mutually inverse and establish  isomorphisms between $V(1)$ and $V(2)$. Clearly, 
$\algA_c=V(1)\oplus V(2)$.  Let $V=V(1)$, for any $v\in V$ we denote  $v(1)=v,\ v(2)=\pi_{12}(v)$, then $\Ca(v)=F\cdot v(1)+F\cdot v(2)\cong \Ca$.

\begin{pro}\label{pr3}
For any $u,v\in V$ we have
\bes
\Ca(u)\cdot \Ca(v) = \la u,v\ra\,{\HH}
\ees
 where $\la ,\ra : V\times V\rightarrow Z(\algA)$ is a skew-symmetric bilinear mapping.
In particular,  $\Ca(v)^2=0$ for any $v\in V$.
\end{pro}
\begin{dem}
 We have by identities of right and left alternativity
\bes
(u(1)v(1))e_{11}&=&-(u(1)e_{11})v(1)+u(1)(v(1)\circ e_{11})=u(1)v(1),\\
e_{11}(u(1)v(1))&=&-u(1)(e_{11}v(1))+(e_{11}\circ u(1))v(1)=-u(1)v(1)+u(1)v(1)=0,
\ees 
which shows that $u(1)v(1)\in e_{22}\algA_ae_{11}=Z_ae_{21}$, hence $u(1)v(1)=ze_{21}$ for some $z\in Z_a$.
 Furthermore, we have
 \bes
 ze_{22}&=&(ze_{21})e_{12}=(u(1)v(1))e_{12}=-(u(1)e_{12})v(1)+u(1)(v(1)\circ e_{12})=u(2)v(1),\\
 ze_{11}&=&e_{12}(ze_{21})=e_{12}(u(1)v(1))=-u(1)(e_{12}v(1))+(e_{12}\circ u(1))v(1)=-u(1)v(2),\\
 ze_{12}&=&e_{12}(ze_{22})=e_{12}(u(2)v(1))=-u(2)(e_{12}v(1))+(e_{12}\circ u(2))v(1)=-u(2)v(2),
 \ees
which proves that $\Ca(u)\cdot\Ca(v)=z\,{\HH}$.

Since $z\in Z_a$, we have $[z,\HH]=[z,\algA_c]=0$. Hence in order to prove that $z\in Z(\algA)$, it remains to show that $[z,Z_a]=0$.
Observe that $z=z(e_{11}+e_{22})=u(2)v(1)-u(1)v(2)\in \algA_c^2$. Therefore,
\[ \ 
[z,Z_a]\subset [\algA_c^2,Z_a]\subseteq \algA_c[\algA_c,Z_a]+[\algA_c,Z_a]\algA_c+3(\algA_c,\algA_c,Z_a)=0.
\] \ 
Finally, denote $z=z(u,v)$ and consider 
$$
e_{11}(u(1)\circ v(1))= (e_{11}u(1)) v(1)+(e_{11}v(1)) u(1)=u(1)\circ v(1).
$$
On the other hand, $e_{11}(u(1)\circ v(1))= e_{11}((z(u,v)+z(v,u))e_{21})=0.$ 
Hence $u(1)\circ v(1)=0$ and $z(u,v)=-z(v,u)$. 
Denote $\la u,v\ra=z(u,v)$, then we have as above
\bee\label{id_bilin}
\la u,v\ra=z(u,v)=u(2)v(1)-u(1)v(2),
\eee
wich proves that $\la u,v\ra$ is a bilinear function of $u,v$.

\end{dem}

\begin{lemma}\label{lem_uvw}
For any $u,v,w,t\in V$ the following identities hold
\bee
\la u,v\ra w+\la v,w\ra u+ \la w,u\ra v=0,\label{e25}\\
 \la u,v\ra \la w,t\ra+\la v,w\ra \la u,t\ra+ \la w,u\ra \la v,t\ra=0.\label{e26}
\eee
\end{lemma}
\begin{dem}  
 Recall that in the proof of Proposition $\ref{pr3}$ we obtained the equalities
\bee
u(1)v(1)&=&\la u,v\ra \,e_{21},\label{uv_11}\\ 
u(1)v(2)&=&-\la u,v\ra \,e_{11},\label{uv_12}\\
u(2)v(1)&=&\la u,v\ra \,e_{22},\label{uv_21}\\ 
u(2)v(2)&=&-\la u,v\ra\, e_{12}\label{uv_22}.
\eee
Therefore, using the fact that $\la u,v\ra\in Z(A)$, by the linearized  right alternative identity  we have
\bes
0&=&(u(1),v(1),w(2))+(u(1),w(2),v(1))\\
&=&\la u,v\ra e_{21}w(2)+ u(1) \la v,w\ra e_{11}-\la u,w\ra e_{11} v(1)-u(1)\la w,v\ra e_{22}\\
&=&\la u,v\ra w+0-\la u,w\ra v-\la w,v\ra u=\la u,v\ra w+\la w,u\ra v+\la v,w\ra u.
\ees
which proves \eqref{e25}.
Multiplying $(\ref{e25})$ by the element $t\in V$, we get  \eqref{e26}.

\end{dem}

\begin{coro}\label{coro_aij}
Let $\{v_i\,|\,i\in I\}$ be a basis of the space $V$ and $u_{ij}=\la v_i,v_j\ra\in Z(\algA)$, then the elements $u_{ij}$ satisfy the Pl\"ucker relations
\bee\label{id_aij}
u_{ij}=-u_{ji},\ \ u_{ij}u_{kl}+u_{ik}u_{lk}+u_{il}u_{jk}=0.
\eee
\end{coro}

An example of a family of elements $u_{ij}=-u_{ji}$  satisfying relations \eqref{id_aij} may be obtained by taking in an associative  commutative algebra $K$ elements $a_1,\ldots,a_n$ and setting $u_{ij}=a_i-a_j$.

\smallskip
Another example, which we will use later, is the coordinate algebra of grassmanian $G_{2,n}$ (see, for example, \cite[vol.1,\,p.42]{Sha}).

\begin{lemma} \label{l8}
	 Consider the algebra of polynomials $F[x_1,\dots,x_n;y_1,\dots,y_n]$, and let 
	$\alpha_{ij}=\mbox{det}\begin{bmatrix}
	x_{i}&y_{i}\\
	x_{j}&y_{j}\\
	\end{bmatrix}\in F[x_1,\dots,x_n;y_1,\dots,y_n]$. Then the elements $\a_{ij}=-\a_{ij}$ satisfy relations \eqref{id_aij}.
	Moreover, the algebra $F[\a_{ij}\,|\,1\leq i<j\leq n]$ is a free algebra modulo  relations \eqref{id_aij}.
\end{lemma}
\begin{dem}
Firstly, one can easily check that the elements $\a_{ij}$ satisfy relations \eqref{id_aij}. 
Furthermore, it follows from the relation
$$
\a_{12}\a_{ij}+\a_{1i}\a_{j2}+\a_{1j}\a_{2i}=0,
$$ 
that $\a_{ij}$ for $i,j>2$  lies in the algebra $F[\a_{1i},\,\a_{2j},\,\a_{12}^{-1}]\subset F(x_1,\ldots, x_n;y_1,\ldots,y_n)$. Therefore, 
\bee
F[\a_{ij}\,|\,1< i\leq n, \,2< j\leq n]\subseteq F[\a_{12},\ldots,\a_{1n};\a_{23},\ldots,\a_{2n},\a_{12}^{-1}]. \label{e27}
\eee
Observe that $y_2=\frac{1}{x_1}\a_{12}+\frac{x_2 y_1}{x_1}$, hence $F(x_1,x_2,y_1,y_2)=F(x_1,x_2,y_1,\a_{12})$.

 Similarly, resolving with respect to $x_n,\,y_n$ the system 
\begin{equation*}
\begin{split}
\a_{1n}& =x_1 y_n-y_1 x_n,\\
\a_{2n}& =x_2 y_n-y_2 x_n
\end{split}
\end{equation*}
	we get 
\bes
	x_n&=&\frac{x_2\a_{1n}-x_1\a_{2n}}{\a_{12}},\\
	y_n&=&\frac{y_2\a_{1n}-y_1\a_{2n}}{\a_{12}}; 
\ees
hence
	\begin{equation*}
	x_n,\, y_n\in F(\a_{1n},\,\a_{2n},\,x_1,\,x_2,\,y_1,\,y_2)=F(\a_{1n},\,\a_{2n},\,x_1,\,x_2,\,y_1,\,\a_{12}).
	\end{equation*}
	Therefore,
	\begin{equation*}
    F(x_1,\ldots,x_n,y_1,\ldots,y_n)=F(x_1,x_2,y_1,\a_{12},\ldots,\a_{1n},\a_{23},\ldots, \a_{2n}).
	\end{equation*}
	and $\mbox{tr.deg}\,F(x_1,x_2,y_1,\a_{12},\ldots,\a_{1n},\a_{23},\ldots, \a_{2n})=2n$,
	 which means that the elements $\a_{12},\ldots,\a_{1n},\a_{23},\ldots, \a_{2n}$ are algebraically independent.

\medskip
Now, let $F[u_{ij}]$ be a free algebra modulo  relations \eqref{id_aij}. Consider
the epimorphism $\pi:F[u_{ij}]\longrightarrow F[\alpha_{ij}]$; $u_{ij}\longmapsto \alpha_{ij}$. We will prove that $\ker\pi=0$. 
Let $f(u_{12},\ldots,u_{(n-1)n})\in \mbox{ker}\,\pi$, that is,  $f(\a_{12},\dots,\a_{(n-1)n})=0$. Inclusions \eqref{e27} follow from relations \eqref{id_aij}, hence they are valid in the algebra $F[u_{ij}]$ as well. Therefore,  there exists $k$ such that 
$$
u_{12}^{k}f(u_{12},\dots,u_{(n-1)n})=g(u_{12},\dots,u_{2n}) 
$$
for some  $g(u_{12},\dots,u_{2n})\in F[u_{12},\dots,u_{2n}]$. Clearly,  
 $g(\a_{12},\dots,\a_{2n})=0$. Since the elements  $\alpha_{12},\ldots,\alpha_{2n}$ are algebraically independent, we have $g=0$.
But the algebra $F[u_{ij}]$ is a domain (see, for example, \cite[Chapter 8]{Muk}),  therefore $f=0$, proving the lemma.
 
\end{dem}

\smallskip

Recall that by Corollary \ref{c2} we have $\algA_a\cong M_2(Z_a)$, hence $\algA=M_2(Z_a)\oplus V(1)\oplus V(2)$.

\begin{pro}\label{pro_product}
Let $X,Y\in \algA,\ X=X_a+x(1)+y(2),\ Y=Y_a+z(1)+t(2)$, where $X_a=\left(\begin{array}{cc}
a&b\\
c&d
\end{array}\right),\ Y_a=\left(\begin{array}{cc}
e&f\\
g&h
\end{array}\right),\ a,b,c,d,e,f,g,h\in Z_a,\ x,y,z,t\in V$. Then the product $XY$ iz given by
\bes
XY=X_aY_a+\left(\begin{array}{cc}
-\la x,t\ra&-\la y,t\ra\\
\la x,z\ra&\la y,z\ra
\end{array}\right)+(az+ct+hx-gy)(1)+(bz+dt-fx+ey)(2).
\ees
\end{pro}
\begin{dem}
The proof follows from identities \eqref{id_cay}, \eqref{uv_11} - \eqref{uv_22}, and Lemma \ref{l4}.

\end{dem}

We can make the formula defining the product in $\algA$ more transparent by using the following notation: for $u,v\in V$ we denote
\bes
(u,v)=u(1)+v(2).
\ees
With this notation, using usual matrix multiplication and the fact that $[Z_a,V_c]=0$, we have for $ X=X_a+(x,y),\ Y=Y_a+(z,t)$
\bee\label{id_product}
XY=X_aY_a+\left(\begin{array}{cc}
-\la x,t\ra&-\la y,t\ra\\
\la x,z\ra&\la y,z\ra
\end{array}\right)+(z,t)X_a+(x,y) (Y_a)^*,
\eee
where ${\left(\begin{array}{cc}
a&b\\
c&d
\end{array}\right)^*}=\left(\begin{array}{cc}
d&-b\\
-c&a
\end{array}\right)$.

In the next section we will prove that Proposition \ref{pro_product} in fact describes all unital alternative extensions of the algebra $M_2(F)$. 

\section{The Main Theorem}


Let $\algB$ be an associative unital algebra and $V$ be a left $\algB$-module such that $[\algB,\algB]$ annihilates $V$. Clearly, in this case $V$ has a structure of a commutative $\algB$-bimodule with $v\cdot b=b\cdot v,\ v\in V,\,b\in \algB$. Assume that there exists a $\algB$-bilinear skew-symmetric mapping $\la,\ra :V^2\rightarrow \algB$ such that $\la V,V\ra\subseteq Z(\algB)$ and formula \eqref{e25} holds for any $u,v,w\in V$.

\smallskip

Let $\algA=M_2(\algB)\oplus V^{2}$, where $V^{2}=\{(u,v)\,|\,u,v\in V\}\cong V\oplus V$. 
Let $X,Y\in \algA,$ $X=X_a+(x,y),\ Y=Y_a+(z,t)$, where $X_a,\ Y_a\in M_2(\algB)$ and $(x,y), (z,t)\in V^2$. Define  a product in $\algA$ by formula \eqref{id_product}:
\bes
XY=X_aY_a+\left(\begin{array}{cc}
-\la x,t\ra&-\la y,t\ra\\
\la x,z\ra&\la y,z\ra
\end{array}\right)+(z,t)X_a+(x,y)(Y_a)^*.
\ees

\begin{theo}\label{th_main}
The algebra $\algA$ with the product defined above is an alternative unital algebra containing $M_2(F)$ with the same unit. Conversely, every unital alternative algebra that contains the matrix algebra $M_2(F)$ with the same unit has this form.
\end{theo}
\begin{dem}
The second part of the theorem follows from Proposition \ref{pro_product} with $B=Z_a$. Let us now prove that $\algA$ is  alternative.

Let us first prove that $V^2$ is a right alternative bimodule over $M_2(\algB)$. Let $A,B\in M_2(\algB)$, $(x,y)\in V^2$. One can easily check that 
$(x,y)({(AB)^*}-B^* A^*)=0$ (since $V[\algB,\algB]=0$). 
Therefore, 
\bes
((x,y),A,A)=((x,y) A^*) A^*-(x,y){(A^2)^*}=0.
\ees
Furthermore,
\bes
&&(A,(x,y),B)+(A,B,(x,y))=\\
&=&((x,y)A) B^*-((x,y) B^*)A+(x,y)(AB)-((x,y) B^*) A^*\\
&=&(x,y)(A B^*- B^*A+AB- B^* A^*)=(x,y)[A,B+ B^*]=0
\ees
since $B+ B^*=tr(B)$ commutes with $A$ on $V^2$. Therefore, $V^2$ is a right alternative bimodule over $M_2(\algB)$.
Now let $A=\left(\begin{array}{cc}
a&b\\
c&d\end{array}\right)$ with $a,b,c,d\in\algB$,
consider 
\bes
(A,(x,y),(x,y))&=&((x,y)A)\cdot (x,y)-\left(\begin{array}{cc}
-\la x,y\ra&0\\
0&\la y,x\ra\end{array}\right)A\\
&=&(xa+yc,xb+yd)\cdot (x,y)+\la x,y\ra A\\
&=&\left(\begin{array}{cc}
-\la xa+yc,y\ra&-\la xb+yd,y\ra\\
\la xa+yc,x\ra&\la xb+yd,x\ra
\end{array}\right)+\la x,y\ra A\\
&=&\left(\begin{array}{cc}
-\la xa,y\ra&-\la xb,y\ra\\
\la yc,x\ra&\la yd,x\ra
\end{array}\right)+\la x,y\ra A=-\la x,y\ra A+\la x,y\ra A=0.
\ees
Furthermore,
\bes
&&((x,y),A,(u,v))+((x,y),(u,v),A)=\\
&=&((x,y) A^*)\cdot (u,v)-(x,y)\cdot ((u,v)A)+((x,y)\cdot (u,v))A-(x,y)\cdot ((u,v) A^*)\\
&=&(xd-yc,-xb+ya)\cdot (u,v)-(x,y)\cdot (ua+vc,ub+vd)\\
&+&\left(\begin{array}{cc}
-\la x,v\ra&-\la y,v\ra\\
\la x,u\ra& \la y,u\ra\end{array}\right)A-(x,y)\cdot(ud-vc,-ub+va)\\
&=&\left(\begin{array}{cc}
-\la xd-yc,v\ra& -\la -xb+ya,v\ra\\
\la xd-yc,u\ra &\la -xb+ya,u\ra\end{array}\right)
-\left(\begin{array}{cc}
-\la x,ub+vd\ra& -\la y,ub+vd\ra\\
\la x,ua+vc\ra &\la y,ua+vc\ra\end{array}\right) \\
&+& \left(\begin{array}{cc}
-\la x,v\ra a-\la y,v\ra c& -\la x,v\ra b-\la y,v\ra d\\
\la x,u\ra a+\la y,u\ra c &\la x,u\ra b+\la y,u\ra d\end{array}\right)
-\left(\begin{array}{cc}
-\la x,-ub+va\ra& -\la y,-ub+va\ra \\
\la x,ud-vc\ra &\la y,ud-vc\ra \end{array}\right)\\
&=&\left(\begin{array}{cc}
X_{11}&X_{12}\\
X_{21}&X_{22}\end{array}\right).
\ees
We have
\bes
X_{11}&=&-\la xd-yc,v\ra+\la x,ub+vd\ra-\la x,v\ra a-\la y,v\ra c+\la x,-ub+va\ra\\
&=&-\la x,v\ra d+\la y,v\ra c+\la x,u\ra b+\la x,v\ra d-\la x,v\ra a-\la y,v\ra c -\la x,u\ra b+\la x,v\ra a=0,
\ees
and similarly
\bes
X_{12}&=&-\la -xb+ya,v\ra+\la y,ub+vd\ra-\la x,v\ra b-\la y,v\ra d+\la y,-ub+va\ra\\
&=&\la x,v\ra b-\la y,v\ra a+\la y,u\ra b+\la y,v\ra d-\la x,v\ra b-\la y,v\ra d -\la y,u\ra b+\la y,v\ra a=0,\\[0,5mm]
X_{21}&=&\la xd-yc,u\ra-\la x,ua+vc\ra+\la x,u\ra a+\la y,u\ra c-\la x,ud-vc\ra \\
&=&\la x,u\ra d-\la y,u\ra c-\la x,u\ra a-\la x,v\ra c+\la x,u\ra a+\la y,u\ra c-\la x,u\ra d+\la x,v\ra c =0,\\[0.5mm]
X_{22}&=&\la -xb+ya,u\ra-\la y,ua+vc\ra+\la x,u\ra b+\la y,u\ra d-\la y,ud-vc\ra \\
&=&\la x,u\ra d-\la y,u\ra c-\la x,u\ra a-\la x,v\ra c+\la x,u\ra a+\la y,u\ra c-\la x,u\ra d+\la x,v\ra c =0.
\ees
Finally,
\bes
&&((x,y),(z,t),(z,t))=\\
&=&\left(\begin{array}{cc}
-\la x,t\ra&-\la y,t\ra\\
\la x,z\ra&\la y,z\ra
\end{array}\right)\cdot (z,t)-(x,y)\cdot 
\left(\begin{array}{cc}
-\la z,t\ra&0\\
0&\la t,z\ra
\end{array}\right) \\
&=&(z,t)\left(\begin{array}{cc}
-\la x,t\ra&-\la y,t\ra\\
\la x,z\ra&\la y,z\ra
\end{array}\right)+\la z,t\ra (x,y)\\
&=&(-\la x,t\ra z+\la x,z\ra t+\la z,t\ra x,-\la y,t\ra z+\la y,z\ra t+\la z,t\ra y)\\
&=&(\la t,x\ra z+\la x,z\ra t+\la z,t\ra x,\la t,y\ra z+\la y,z\ra t+\la z,t\ra y)\stackrel{\eqref{e25}}=(0,0).
\ees
Therefore, $\algA$ is right alternative. Similarly, one can prove that $\algA$ is left alternative.
 
\end{dem}

\section{Examples}

\subsection{Algebra of octonions.}

Let $\algB$ be a unital associative  commutative algebra and $\algA=\OO(\algB)$ be a split octonion algebra over $\algB$.
In this case, $\algA=M_2(\algB)\oplus vM_2(\algB)$ with $v^2=1,\ \algA_a=M_2(\algB),\ \algA_c=vM_2(\algB)$, $Z_a=Z(\algA)=\algB$.
Take $V=\algB^2=\{(a,b)\,|\,a,b\in\algB\}$, $(a,b)(1)=v\left(\begin{array}{cc}
0&0\\ 
-b&a\end{array}\right),\ 
(a,b)(2)=v\left(\begin{array}{cc}
b&-a\\ 
0&0\end{array}\right),$ then we have $\algA=M_2(\algB)\oplus V(1)\oplus V(2)$,  with $\la(a,b),(c,d)\ra=-\det \left(\begin{array}{cc}
a&b\\ 
c&d\end{array}\right)$.
In fact, by \eqref{id_bilin},
\bes
\la (a,b),(c,d)\ra&=&(a,b)(2)\cdot (c,d)(1)-(a,b)(1)\cdot (c,d)(2)\\
&=&v\left(\begin{array}{cc}
b&-a\\ 
0&0\end{array}\right)\cdot v\left(\begin{array}{cc}
0&0\\ 
-d&c\end{array}\right)-v\left(\begin{array}{cc}
0&0\\ 
-b&a\end{array}\right)\cdot v\left(\begin{array}{cc}
d&-c\\ 
0&0\end{array}\right)\\
&=&\left(\begin{array}{cc}
0&0\\ 
-d&c\end{array}\right)\cdot \left(\begin{array}{cc}
0&a\\ 
0&b\end{array}\right)-\left(\begin{array}{cc}
d&-c\\ 
0&0\end{array}\right)\cdot \left(\begin{array}{cc}
a&0\\ 
b&0\end{array}\right)=-\det \left(\begin{array}{cc}
a&b\\ 
c&d\end{array}\right).
\ees
Now for any $u=(a,b),\,v=(c,d),\, w=(e,f)\in V$ we have
\bes
&&\la u,v\ra w+\la v,w\ra u+\la w,u\ra v=\\
&=&-\det \left(\begin{array}{cc}
a&b\\ 
c&d\end{array}\right) (e,f)-\det\left(\begin{array}{cc}
c&d\\ 
e&f\end{array}\right) (a,b)-\det\left(\begin{array}{cc}
e&f\\ 
a&b\end{array}\right) (c,d)\\
&=&-(ad-bc)(e,f)-(cf-de)(a,b)-(eb-fa)(c,d)=(0,0),
\ees
hence $\OO(\algB)$ satisfies \eqref{e25}.

\smallskip

The following Proposition gives conditions under which the algebra $\algA$ from Theorem \ref{th_main} is isomorphic to an octonion algebra $\OO(\algB)$.
\begin{pro}
The unital algebra $\algA=M_2(\algB)\oplus V^2$ from Theorem \ref{th_main} is isomorphic to an octonion algebra $\OO(\algB)$ if and only if there exist $x,y\in V$ such that $\la x,y\ra=1$.
\end{pro}
\begin{dem}
We have already checked that  the algebra $\OO(\algB)$ has form $M_2(\algB)\oplus V^2$, it suffices to note that $\la x,y\ra=1$ for $x=(1,0),\ y=(0,-1)\in V$.

Let now $\algA=M_2(\algB)\oplus V^2$ be such that there exist $x,y\in V$ with $\la x,y\ra=1.$ 
Observe first that for any $u,v\in V,\, a,b\in \algB$ the following equality holds 
\bee\label{[a,b]uv}
[a,b]\la u,v\ra=0.
\eee
In fact, we have
\bes
ab\la u,v\ra&=&a\la bu, v\ra=\la bu,av\ra=b\la u,av\ra=ba\la u,v\ra.
\ees
For any $a,b\in\algB$ we now have $0=[a,b]\la x,y\ra=[a,b],$ hence $\algB$ is commutative. Consider $\algC=M_2(F)+ \Ca(x)+\Ca(y)$.
It follows from Proposition \ref{pr3} and its proof that $\algC$ is a subalgebra of $\algA$ isomorphic to the split octonion algebra $\OO(F)$.
Therefore, by Kaplansky-Jacobson Theorem, $\algA\cong\OO(A)$ for some commutative associative algebra $A$. It follows from \eqref{e25} that $V=\algB\cdot\Ca(x)+\algB\cdot \Ca(y)$ and  $A=\algB$.

\end{dem}

\subsection{Algebras obtained by (commutative) Cayley-Dickson process}

Note that if the mapping $\la ,\ra :V^2\rightarrow Z(\algB)$ is trivial then the algebra $\algA$ is just a split null extension 
of the algebra $M_2(\algB)$ by a bimodule $V^2$. In this case, $V$ may be an arbitrary associative  $\algB$-module, (annihilated by $[\algB,\algB]$ if $\algB$ is not commutative). For instance, when $\algB=F$ and $V=F$ we get in this way the algebra $\algA=M_2(F)\oplus \Ca$.

\smallskip

If the mapping $\la ,\ra :V^2\rightarrow Z(\algB)$ is not trivial then by \eqref{e25} the rank of $V$ as a $\algB$-module is less than 3. 
Observe that the left side of \eqref{e25} is $\algB$-multilinear and skew-symmetric on $u,v,w$. Therefore, it holds when $\Lambda^3(V_{\algB})=0$. In particular it holds if the rank  of $V$ is less or equal to 2. 
If $V\subseteq \algB\cdot x$ then the mapping $\la,\ra$ is trivial by skew-symmetry. Let us consider now the case when $V$ is a 2-generated $\algB$-module.

\smallskip

Let $A$ be an associative commutative algebra and $\alpha\in A$. Denote by $CD(M_2(A),\alpha)$ the algebra $M_2(A)\oplus vM_2(A)$ with a product defined by the following analogue of \eqref{CD_product}: 
\bee\label{CD_product1}
{{a\cdot b=ab,\ a\cdot vb=v(a^*b),\ vb\cdot a=v(ab), \ va\cdot vb=\alpha (b a^*),}}
\eee
where $a,b\in M_2(A),\  a\mapsto  a^*$ is the symplectic involution in $M_2(A)$. The algebra $CD(M_2(A),\alpha)$ is an alternative algebra containing $M_2(A)$ with the same unit. We will call it the {\em algebra obtained from $M_2(A)$ by the Cayley-Dickson process with a parameter $\alpha$}. The algebra $CD(M_2(A),\a)$ is an octonion algebra if and only if the parameter $\a$ is invertible in $A$.

\begin{theo}\label{th_CD}
Let $\algB$ be a unital commutative algebra, $V=\algB^2$ and $\la ,\ra:V^2\rightarrow \algB$ be a skew-symmetric $\algB$-bilinear mapping. Then the algebra $\algA=M_2(\algB)\oplus V^2$ is isomorphic to an algebra $CD(M_2(\algB),\alpha)$ where $\alpha=-\la (1,0),(0,1)\ra$. Conversely, every algebra $CD(M_2(A),\alpha)$ has this form. 
\end{theo}
\begin{dem}
Let $\algA=CD(M_2(A),\alpha)$. Take $V=A^2=\{(a,b)\,|\,a,b\in A\}$,  $(a,b)(1)=v\left(\begin{array}{cc}
0&0\\ 
-b&a\end{array}\right)$,\ 
$(a,b)(2)=v\left(\begin{array}{cc}
b&-a\\ 
0&0\end{array}\right)\in vM_2(A)$, then we have, as before, $\algA=M_2(A)\oplus V(1)\oplus V(2)$, 
  with $\la(a,b),(c,d)\ra=-\a\det \left(\begin{array}{cc}
a&b\\ 
c&d\end{array}\right)$. In particular, $\la (1,0),(0,1)\ra=-\a$.

Conversely, let $\algA=M_2(\algB)\oplus V^2$ where $V\cong\algB^2$ and $\la(1,0),(0,1)\ra=-\a$. 
Define the mapping $\f:V^2=V(1)\oplus V(2)\rightarrow vM_2(\algB)\subset CD(M_2(\algB),\a)$ by sending, for any $a,b\in \algB$
\bes
(a,b)(1)\mapsto v\left(\begin{array}{cc} 0&0\\-b&a\end{array}\right),\ 
(a,b)(2)\mapsto v\left(\begin{array}{cc} b&-a\\0&0\end{array}\right).
\ees
It is easy to see that $\f$ is an isomorphism of alternative $M_2(\algB)$-bimodules. 
Furthermore, let $x=(a,b),\, y=(c,d)\in V=\algB^2$, then we have 
\bes
\la x,y\ra&=&\la (a,b),(c,d)\ra=\la a(1,0)+b(0,1),c(1,0)+d(0,1)\ra\\
&=&(ad-bc)\la(1,0),(0,1)\ra=-\a(ad-bc).
\ees
Let $z=(e,f),\,t=(g,h)\in V$, then we have by \eqref{id_product}
\bes
(x,y)(z,t)= \left(\begin{array}{cc}
-\la x,t\ra&-\la y,t\ra\\
\la x,z\ra&\la y,z\ra
\end{array}\right)=-\a\left(\begin{array}{cc}
-ah+bg&-ch+dg\\
af-be&cf-de
\end{array}\right).
\ees
On the other hand,
\bes
\f(x,y)\cdot \f(z,t)&=&v\left(\begin{array}{cc} d&-c\\-b&a\end{array}\right)\cdot v\left(\begin{array}{cc} h&-g\\-f&e\end{array}\right)=\a\left(\begin{array}{cc} h&-g\\-f&e\end{array}\right)\cdot
\left(\begin{array}{cc} a&c\\b&d\end{array}\right)\\
&=&\a\left(\begin{array}{cc} ah-bg&ch-dg\\-af+be&-cf+ed\end{array}\right).
\ees
Therefore, the mapping $id+\f:\algA= M_2(\algB)\oplus V^2\rightarrow CD(M_2(\algB),\a)=M_2(\algB)\oplus vM_2(\algB)$ is an isomorfism.

\end{dem}

\subsection{Algebras obtained by noncommutative Cayley-Dickson process}

Let us now generalize the Cayley-Dickson process for non-commutative coefficient algebras.
Let $A$ be a unital associative algebra, not necessarily commutative, $\a\in A$ such that $\a A\subseteq Z(A)$.
Denote $NCD(M_2(A),\a)=M_2(A)\oplus vM_2(\bar A)$, where $\bar A=A/[A,A]A$, and define a product in it by setting
\bee\label{NCD_product}
{{a\cdot b=ab,\ a\cdot v\bar b=v(\bar a^*\bar b),\ v\bar b\cdot a=v(\bar a\bar b), \ v\bar a \cdot v\bar b=\alpha (b_1a_1^*),}}
\eee
where $a,b\in M_2(A),\ \bar a,\bar b$ are their images in $M_2(\bar A)$, $\bar a\mapsto \bar a^*$ is the symplectic involution in $M_2(\bar A)$, and $a_1^*,b_1\in M_2(A)$  are some pre-images of $\bar a^*,\bar b$ under the canonical epimorphism $M_2(A)\rightarrow M_2(\bar A)$. Observe that the last product in \eqref{NCD_product} is defined correctly since $\a[A,A]=0$.

\begin{pro}\label{pro_NCD}
The algebra $NCD(M_2(A),\a)$ is a unital alternative algebra that contains $M_2(A)$ with the same unit.
\end{pro}
\begin{dem}
Denote $I=[A,A]A$, then $M_2(I)$ is an ideal of $NCD(M_2(A),\a)$ which annihilates $vM_2(\bar A)$ and is annihilated by $\a$; moreover, $NCD(M_2(A),\a)/I\cong CD(M_2(\bar A),\a)$.
Therefore, the $M_2(A)$-bimodule $vM_2(\bar A)$ is in fact an $M_2(\bar A)$-bimodule, and since the algebra $CD(M_2(\bar A),\a)$ is alternative,  $vM_2(\bar A)$ is an alternative $M_2(A)$-bimodule. 
In this way, it suffices to check the alternativity identities only when we have at least two arguments belonging to $vM_2(\bar A)$. 

For any $a,b\in M_2(A)$ we have
\bes
(a,v\bar b,v\bar b)&=&(v(\bar a^*\bar b))\cdot v\bar b-a(\a bb_1^*)=\a\, b(b_1^*a)-\a\, a(bb_1^*),
\ees
where $\overline{b_1^*}=\bar b^*$. Consider 
\bes
\overline{b(b_1^*a)- a(bb_1^*)}=\bar b(\bar b^*\bar a)-\bar a(\bar b\bar b^*)=(\det\bar b)\bar a-\bar a(\det\bar b)=\bar 0.
\ees
Thus $b(b_1^*a)- a(bb_1^*)\in M_2(I)$ and $\a(b(b_1^*a)- a(bb_1^*))=0$.

Furthemore,
\bes
(v\bar a,v\bar b,v\bar b)&=&(\a ba_1^*)\cdot v\bar b-v\bar a\cdot (\a bb_1^*)=\a v((\bar a\bar b^*)\bar b-(\bar b\bar b^*)\bar a)=0.
\ees
Finally, consider, for $c\in M_2(A)$,
\bes
(v\bar a,v\bar b,c)+(v\bar a,c,v\bar b)&=&\a\,ba_1^*\cdot c-v\bar a\cdot v(\bar c\bar b)+v(\bar c\bar a)\cdot v\bar b-v\bar a\cdot v(\bar c^*\bar b)\\
&=&\a\,(ba_1^*\cdot c- cb\cdot a_1^*+b\cdot a_1^*c_1^*-c_1^*b\cdot a_1^*)
\ees
We have
\bes
\overline{ba_1^*\cdot c- cb\cdot a_1^*+b\cdot a_1^*c_1^*-c_1^*b\cdot a_1^*}&=&\bar b\bar a^*\cdot \bar c-\bar c\bar b\cdot \bar a^*+\bar b\cdot \bar a^*\bar c^*-\bar c^*\bar b\cdot \bar a^*\\
&=&\bar b\bar a^*\,t(\bar c)-t(\bar c)\,\bar b\bar a^*=\bar 0.
\ees
Hence $ba_1^*\cdot c- cb\cdot a_1^*+b\cdot a_1^*c_1^*-c_1^*b\cdot a_1^*\in M_2(I)$ and $\a\,(ba_1^*\cdot c- cb\cdot a_1^*+b\cdot a_1^*c_1^*-c_1^*b\cdot a_1^*)=0$.

\smallskip
We have proved that the algebra $NCD(M_2(A),\a)$ is right alternative. Similarly, one can check that it is left alternative.

\end{dem}

 Now we can generalize Theorem \ref{th_CD} to the case when $\algB$ is not commutative.
\begin{theo}\label{th_NCD}
Let $\algB$ be a unital  associative algebra, $\overline{\algB}=\algB/[\algB,\algB]\algB$, $V=\overline{\algB}^2$ and $\la ,\ra:V^2\rightarrow \algB$ be a skew-symmetric $\algB$-bilinear mapping. Then the algebra $\algA=M_2(\algB)\oplus V^2$ is isomorphic to an algebra $NCD(M_2(\algB),\alpha)$ where $\alpha=-\la (1,0),(0,1)\ra$. Conversely, every algebra $NCD(M_2(A),\alpha)$ has this form. 
\end{theo}
\begin{dem}
Let first $\algA=NCD(M_2(A),\a)$. Denote $\bar A=A/[A,A]A$ and take $V = \bar A^2 = {(\bar a,\bar b)\,|\,a,b\in A},\ (\bar a,\bar b)(1) = v\left(\begin{array}{cc}
0&0\\ 
-\bar b&\bar a\end{array}\right)$,\ 
$(\bar a,\bar b)(2)=v\left(\begin{array}{cc}
\bar b&-\bar a\\ 
0&0\end{array}\right)\in vM_2(\bar A)$; then we have, as before, $\algA=M_2(A)\oplus V(1)\oplus V(2)$, 
  with $\la(\bar a,\bar b),(\bar c,\bar d)\ra=-\a(ad-bc)$. In particular, $\la (\bar 1,\bar 0),(\bar 0,\bar 1)\ra=-\a$.
\smallskip

Conversely, let $\algA=M_2(\algB)\oplus V^2$ where $V\cong\overline{\algB}^2,\ \overline{\algB}=\algB/[\algB,\algB]\algB$ and $\la(\bar 1,\bar 0),(\bar 0,\bar 1)\ra=-\a$. 
Define the mapping $\f:V^2=V(1)\oplus V(2)\rightarrow vM_2(\overline{\algB})\subset NCD(M_2(\algB),\a)$ by sending, for any $a,b\in \algB$
\bes
(\bar a,\bar b)(1)\mapsto v\left(\begin{array}{cc} 0&0\\-\bar b&\bar a\end{array}\right),\ 
(\bar a,\bar b)(2)\mapsto v\left(\begin{array}{cc} \bar b&-\bar a\\0&0\end{array}\right).
\ees
Then, as in the proof of Theorem \ref{th_CD}, one can easily see that the mapping
\bes
id+\f:\algA= M_2(\algB)\oplus V^2\rightarrow NCD(M_2(\algB),\a)=M_2(\algB)\oplus vM_2(\overline{\algB})
\ees
 is an isomorfism. 

\end{dem}

Algebras of type $CD(M_2(A),\a)$ can be constructed for any commutative algebra $A$ and any $\a\in A$. 
For algebras of type $NCD(M_2(A),\a)$ one have to check the condition $[\a A,A]=0$. For instance,
one can take $A=F\la x,y\,|\,y[x,y]=[x,y]y=0\ra$, with $\a=y^2.$

\subsection{The case when $V$ is not 2-generated}

In all the examples considered above the $\algB$-module $V$ was 2-generated. 
Here we will give an example when $V$ is 3-generated. 

Let $\algB$ be a commutative unital algebra, $a,b,c\in\algB,\ V=\algB^3/I$ where $I=\algB\cdot(a,b,c)$; denote $e_1=(1,0,0)+I,\,e_2=(0,1,0)+I,\, e_3=(0,0,1)+I$.
Then we have  $V=\algB\cdot e_1+\algB\cdot e_2+\algB\cdot e_3$ where $a\cdot e_1+b\cdot e_2+c\cdot e_3=0$.
 Define a $\algB$-bilinear skew-symmetric mapping $\la,\ra:V\times V\rightarrow \algB$ by setting
\bes
\la e_1,e_2\ra=c,\ \la e_2,e_3\ra=a,\ \la e_3,e_1\ra=b.
\ees
 One can easily check that the mapping $\la,\ra$ is defined correctly. Moreover, we have
\bes
\la e_1,e_2\ra e_3+\la e_2,e_3\ra e_1+\la e_3,e_1\ra e_2=c\cdot e_3+a\cdot e_1+b\cdot e_2=0,
\ees 
that is, identity \eqref{e25} is true for $u=e_1,\,v=e_2,w=e_3$. Since the left side of \eqref{e25} is skew-symmetric and multilinear on $u,v,w$, it follows that \eqref{e25} is valid in $V$.
By Theorem \ref{th_main}, the algebra $\algA=M_2(\algB)\oplus V^2$  is a unital alternative algebra containing $M_2(\algB)$ as a unital subalgebra.

\smallskip

Observe that taking here $a=b=0$, we will get the algebra $CD(\algB,c)$ from the Theorem \ref{th_CD}. 

Moreover, following the scheme from the previous section, this construction can be extended for noncommutative algebra $\algB$. One has only to choose the elements $a,b,c\in\algB$ such that 
 $a\algB+b\algB+c\algB\subset Z(\algB)$.

\section{Open questions}

{\bf 1.} The first natural question which we left open is the case when the algebra $\HH$ is not split, that is, when $\HH$ is a division algebra.
This case is more complicated since while $\Ca_i\cdot \Ca_j=\a_{ij}\HH$, the product $\Reg_i\HH\cdot \Reg_j\HH=\sum_{k=1}^4\a_{ij}^k\HH$ for some $\a_{ij}^k\in Z(\algA)$, and instead of Pl\"ucker relations \eqref{id_aij} the elements $\a_{ij}^k$ satisfy more complicated system of relations. We plan to consider this case in a forthcoming paper.

\medskip

{\bf 2.} An interesting question is to study the alternative algebras that contain $\HH$ (or {\em $\HH$-algebras}) from categorical point of view. Clearly, the class of $\HH$-algebras form a category, with morphisms being the homomorphisms acting identically on $\HH$. Given an $\HH$-bimodule $V$, the free $\HH$-algebra over $V$ or {\em tensor algebra $\HH[V]$ of the bimodule $V$} plays a role of a free object in this category. When $V=V_a$ is associative, $V=\oplus_{i=1}^m\Reg_i\HH$ and  the algebra $\HH[V]$ is associative and is isomorphic to $\HH\otimes F\la x_1,\ldots,x_m\ra$ where $F\la x_1,\ldots,x_m\ra$ is the free associative algebra on $m$ generators.
 
When $V=V_c$ is a Cayley $\HH$-bimodule, $V=\oplus_{i=1}^m\Ca_i$, the situation is not so clear even in the split case.  For $m=1$, $\HH[V]=\HH\oplus\Ca$ with $\Ca^2=0$ is just a well known  6-dimensional subalgebra of a split Cayley-Dickson algebra; for $n=2$, $\HH[V]\cong CD(M_2(F[\a_{12}]),\a_{12})$, but for $n\geq 3$ the structure of the algebra $\HH[V]$ is not known.

The situation is even more complicated for the mixed case, when $V=V_a\oplus V_c$ with $V_a,\,V_c\neq 0$, again except some trivial cases.   

\section{Acknowledgements}

The paper is a part of the PhD-thesis of the first author, realized at the University of S\~ao Paulo with the support of the CAPES. The second author was partially supported by the CNPq grant  304313/2019-0 and by the FAPESP grant  2018/23690-6.

\end{document}